\magnification=\magstep1
\documentstyle{amsppt}

\topmatter
\title Subadjunction of log canonical divisors II
\endtitle
\author Yujiro Kawamata
\endauthor

\rightheadtext{adjunction}

\address Department of Mathematical Sciences, University of Tokyo, 
Komaba, Meguro, Tokyo, 153, Japan \endaddress
\email kawamata\@ms.u-tokyo.ac.jp\endemail

\abstract
We extend a subadjunction formula of log canonical divisors as in [K3] to 
the case when the codimension of the minimal center is arbitrary
by using the positivity of the Hodge bundles.
\endabstract

\endtopmatter

\document

\head 1. Main Result
\endhead

The following is the main result:

\proclaim{Theorem 1}
Let $X$ be a normal projective variety.
Let $D^o$ and $D$ be effective $\Bbb Q$-divisors on $X$ such that
$D^o < D$, $(X, D^o)$ is log terminal, and $(X, D)$ is log canonical.
Let $W$ be a minimal center of log canonical singularities for $(X, D)$.
Let $H$ be an ample Cartier divisor on $X$, 
and $\epsilon$ a positive rational number.
Then there exists an effective $\Bbb Q$-divisor $D_W$ on $W$
such that 
$$
(K_X + D + \epsilon H) \vert_W \sim_{\Bbb Q} K_W + D_W
$$
and that the pair $(W, D_W)$ is log terminal.
In particular, $W$ has only rational singularities.
\endproclaim

We recall the terminology.
A pair $(X, D)$ of a normal variety and
an effective $\Bbb Q$-divisor is said to be {\it log terminal (KLT)} 
(resp. {\it log canonical (LC)}) 
if the following conditions are satisfied:

(1) $K_X + D$ is a $\Bbb Q$-Cartier divisor.

(2) There exists a projective birational morphism $\mu: Y \to X$ from a 
smooth variety $Y$ with a normal crossing divisor $\sum_j E_j$ such that 
a formula
$$
K_Y + \sum_j e_jE_j \sim_{\Bbb Q} \mu^*(K_X + D)
$$
holds with $e_j < 1$ (resp. $\le 1$) for all $j$,
where $\sim_{\Bbb Q}$ denotes the $\Bbb Q$-linear equivalence.

If $e_j = 1$, then the subvariety $\mu(E_j)$ of $X$ is called 
a {\it center of log canonical singularities},
and the discrete valuation of $\Bbb C(X)$ corresponding to the 
prime divisor $E_j$ is called a {\it place of log canonical singularities}.

If $D^o$ and $D$ are effective $\Bbb Q$-divisors on a normal variety $X$ 
such that $D^o < D$, $(X, D^o)$ is KLT, and that $(X, D)$ is LC, then
there exists a minimal element among the centers of log canonical 
singularities for $(X, D)$ with respect to the inclusions (cf. [K2, \S 1]).
Though the argument in [K2, \S 1] on the {\it minimal center of log canonical
singularities} treats only the case where $D^o = 0$, 
it can be easily extended to our case.
In particular, a minimal center of log canonical singularities is always 
normal.

\head 2. Positivity
\endhead

The following positivity result is the key to the proof of Theorem 1.

\proclaim{Theorem 2}
Let $f: X \to B$ be a surjective morphism of smooth projective varieties with
connected fibers.
Let $P = \sum P_j$ and $Q = \sum_{\ell} Q_{\ell}$ be
normal crossing divisors on $X$ and $B$, respectively, such that
$f^{-1}(Q) \subset P$ and 
$f$ is smooth over $B \setminus Q$.
Let $D = \sum_j d_jP_j$ be a $\Bbb Q$-divisor on $X$,
where $d_j$ may be positive, zero or negative, 
which satisfies the following conditions:

(1) $D = D^h + D^v$ such that 
$f: \text{Supp}(D^h) \to B$ is surjective and smooth 
over $B \setminus Q$, and
$f(\text{Supp}(D^v)) \subset Q$.  An irreducible 
component of $D^h$ (resp. $D^v$) is called 
{\it horizontal} (resp. {\it vertical}).

(2) $d_j < 1$ for all $j$.

(3) The natural homomorphism $\Cal O_B \to 
f_*\Cal O_X(\ulcorner - D \urcorner)$ is surjective at the generic point 
of $B$.

(4) $K_X + D \sim_{\Bbb Q} f^*(K_B + L)$ for some 
$\Bbb Q$-divisor $L$ on B.

\noindent
Let 
$$
\align
f^*Q_{\ell} &= \sum_j w_{\ell j}P_j \\
\bar d_j &= \frac {d_j + w_{\ell j} - 1}{w_{\ell j}} 
\text{ if } f(P_j) = Q_{\ell} \\  
\delta_{\ell} &= \text{max }\{ \bar d_j; f(P_j) = Q_{\ell}\} \\
\Delta &= \sum_{\ell} \delta_{\ell}Q_{\ell} \\
M &= L - \Delta.
\endalign
$$
Then $M$ is nef.
\endproclaim

\demo{Proof}
By replacing $D$ by $D - f^*\Delta$, we may assume that $\Delta = 0$.
Then we have an inequality $d_j \le 1 - w_{\ell j}$ for $f(P_j) = Q_{\ell}$, 
and the equality holds for some $j$ for each $\ell$.

By the stable reduction theorem [KKMS] and the covering trick [K1, Theorem 17],
we obtain a semistable reduction in codimension $1$ in the following sense:
there exists a finite morphism $h: B' \to B$ from a smooth projective 
variety $B'$ with a normal crossing divisor 
$Q' = \text{Supp}(h^*Q) = \sum_{\ell'} Q'_{\ell'}$ 
such that the induced morphism $f': X' \to B'$ from a desingularization $X'$ 
of $X \times_B B'$ is semistable over the generic points of $Q'$.
Let $g: X' \to X$ be the induced morphism.
We may assume that $P' = \text{Supp}(g^*P) = \sum_{j'} P'_{j'}$ 
is a normal crossing divisor again:
$$
\CD
X     @<g<< X' \\
@VfVV       @VV{f'}V \\
B     @<h<< B'.
\endCD
$$

Let $Z \subset B'$ be a closed subset of codimension $2$ or larger which is 
contained in $Q'$ and such that $Q' \setminus Z$ is smooth and
$f'$ is semistable over $B' \setminus Z$.
We can define naturally a $\Bbb Q$-divisor $D' = \sum_{j'}d'_{j'}P'_{j'}$ on 
$X'$ such that $K_{X'} + D' \sim_{\Bbb Q} f^{\prime *}(K_{B'} + h^*L)$.
We calculate the coefficients $d'_{j'}$.
If $P'_{j'}$ is horizontal and $g(P'_{j'}) = P_j$, then $d'_{j'} = d_j$.
We are not concerned with those $P'_{j'}$ such that $f'(P'_{j'}) \subset Z$.

We consider the case where $f'(P'_{j'}) = Q'_{\ell'}$.
First, we assume that $g(P'_{j'}) = P_j$ and
$h(Q'_{\ell'}) = Q_{\ell}$.
Let $e_{j'}$ and $e_{\ell'}$ be the ramification indices of $g$ and $h$ at
the generic points of $P'_{j'}$ and $Q'_{\ell'}$, respectively.
Since $f'$ is semistable in codimension $1$, we have
$e_{\ell'} = e_{j'}w_{\ell j}$.
Thus 
$$
d'_{j'} = e_{j'}d_j - (e_{j'} - 1) + e_{\ell'} - 1
= e_{j'}(d_j + w_{\ell j} - 1).
$$
Therefore, we have an inequality $d'_{j'} \le 0$, 
and the equality holds for some $j'$ for each $\ell'$.
The inequality $d'_{j'} \le 0$ holds for any $j'$ such that $f'(P'_{j'}) = 
Q'_{\ell'}$, 
because it is stable under the blow-ups of $X$ (cf. [K3, Remark 6]).

Let $m$ be a positive integer such that $mL$ is a divisor, 
$X'_m = X' \times_{B'} \cdots \times_{B'} X'$ the  
$m$-tuple fiber product of $X'$ over $B'$, and
$f'_m: X'_m \to B'$ the projection.
$X'_m$ has only Gorenstein toric singularities over $B' \setminus Z$.
Let $D'_m = \sum_{i = 1}^m p_i^*D'$ be the sum 
over the projections $p_i: X'_m \to X'$.
Then we have $K_{X'_m} + D'_m \sim_{\Bbb Q} f_m^{\prime *}(K_{B'} + mh^*L)$ 
over $B' \setminus Z$ (cf. [V, Lemma 3.5]).

Let $r$ be the smallest positive integer such that 
$r(K_{X'_m} + D'_m - f_m^{\prime *}(K_{B'} + mh^*L)) \sim 0$ 
over $B' \setminus Z$,
and $\theta$ a corresponding non-zero rational section of 
$\Cal O_{X'_m}(r(K_{X'_m/B'} - mf_m^{\prime *}h^*L))$.
Let $\pi: \tilde X \to X'_m$ be the normalization of the main irreducible 
component $X^{\prime o}_m$ of $X'_m$ in the field 
$\Bbb C(X^{\prime o}_m)(\root{r}\of{\theta})$, and 
$\tilde f: \tilde X \to B'$ the induced morphism.
$\pi$ may ramify along the support of $D'_m$ over $B' \setminus Z$,
and $\tilde f$ may have non-connected fibers.
$\root{r}\of{\theta}$ induces a rational section $\tilde \theta$ of 
$\Cal O_{\tilde X}(K_{\tilde X/B'})$ which is regular along the generic fiber
of $\tilde f$, because the coefficients of $D'$ are less than $1$.

By [K1, Theorem 17] again, 
we construct a finite morphism $h': B'' \to B'$ from a 
smooth projective variety $B''$
with a normal crossing divisor 
$Q'' = \text{Supp}(h^{\prime *}Q') = \sum_{\ell''} Q''_{\ell''}$ 
such that the induced morphism $\tilde f': \tilde X' \to B''$ 
for a desingularization $\mu: \tilde X' \to \tilde X \times_{B'} B''$ 
has unipotent local monodromies around the irreducible components of $Q''$.
Let $\tilde g: \tilde X' \to \tilde X$ be the induced morphism.
By [K1, Theorem 5], $\tilde f'_*\Cal O_{\tilde X'}(K_{\tilde X'/B''})$ is a 
semi-positive locally free sheaf.

Let $X''$ be the normalization of the main irreducible components of the
fiber product $X'_m \times_{B'} B''$.
Let $f'': X'' \to B''$ and $g': X'' \to X'_m$ be the induced morphisms.
Since $f'$ is semistable over $B' \setminus Z$, we have
$g^{\prime *}K_{X'_m/B'} = K_{X''/B''}$ over $B'' \setminus Z'$ for
$Z' = h^{\prime -1}(Z)$.
Thus $\theta' = g^{\prime *}\theta$ is a rational section of 
$\Cal O_{X''}(r(K_{X''/B''} - mf^{\prime \prime *}h^{\prime *}h^*L))$,
and $\tilde X \times_{B'} B''$ coincides with the normalization of $X''$
in the extension $\Bbb C(X'')(\root{r}\of{\theta'})$.
Let $\pi': \tilde X' \to X''$ and 
$\tilde f': \tilde X' \to B''$ be the induced morphisms:
$$
\CD
\tilde X   @<{\tilde g}<< \tilde X' \\
@V{\pi}VV                 @VV{\pi'}V \\
X'_m       @<{g'}<<       X'' \\
@V{f'_m}VV                @VV{f''}V \\
B'         @<{h'}<<       B''.
\endCD
$$
Since $\pi'$ ramifies only along the support of $D'' = g^{\prime *}D'_m$ over 
$B'' \setminus Z'$,
$\tilde X \times_{B'} B''$ has only toric singularities there.
Hence $\mu_*\Cal O_{\tilde X'}(K_{\tilde X'}) 
= \Cal O_{\tilde X \times_{B'} B''}(K_{\tilde X \times_{B'} B''})$ 
over $B'' \setminus Z'$.

$\tilde \theta' = \tilde g^*\tilde \theta$ is a rational section of
$\Cal O_{\tilde X'}(K_{\tilde X'/B''})$.
We may assume that the Galois group $G \cong \Bbb Z_r$ of the covering
$\pi': \tilde X' \to X''$ acts on $\tilde X'$ 
such that a generator of $G$ acts as 
$\tilde \theta' \mapsto \zeta \tilde \theta'$ for a primitive root 
of unity $\zeta$.
The eigen-subsheaf $\Cal L$ of 
$\tilde f'_*\Cal O_{\tilde X'}(K_{\tilde X'/B''})$ corresponding to the
eigenvalue $\zeta$ under the induced action of $G$ is a direct summand, 
hence a semi-positive locally free
sheaf itself.
Since $\tilde \theta'$ gives a rational section of $\Cal L$, 
it is an invertible sheaf by the conditions (2) and (3).  
Moreover, since we have $d'_{j'} \le 0$ with
the equality for some $j'$ for each $\ell'$, 
we have $\Cal L \cong \Cal O_{B''}(mh^{\prime *}h^*L)$.  
Therefore, $L$ is nef.
\hfill $\square$
\enddemo

\demo{Proof of Theorem 1}
There exists an effective $\Bbb Q$-divisor $D'$ which passes through $W$ and
satisfies the following conditions:
$(X, (1-\alpha)D + D')$ is LC for a rational number $\alpha$ such that
$0 < \alpha \ll 1$, $W$ is a minimal center of log canonical singularities
for $(X, (1-\alpha)D + D')$, and there exists only one place of 
log canonical singularities for $(X, (1-\alpha)D + D')$ above $W$.

Let $D_t = (1-\alpha t)D + tD'$ for $0 \le t \le 1$, 
and let $\mu: Y \to X$ be an embedded resolution for the pairs $(X, D_0)$
and $(X, D_1)$ simultaneously. We write
$$
K_Y + E + F_t \sim_{\Bbb Q} \mu^*(K_X + D_t)
$$
where $E$ is the only place of log canonical singularities for 
$(X, D_t)$ above $W$ if $t \ne 0$.
By construction, the coefficients of $F_t \vert_E$ are less than $1$
if $t \ne 0$. Moreover, even if $t = 0$, the coefficients are less than $1$
for vertical components of $F_t \vert_E$ with respect to $\mu$,
because $W$ is a minimal center.

We have $R^1\mu_*\Cal O_Y(- E + \ulcorner - F_t \urcorner) = 0$
by Kawamata-Viehweg vanishing theorem, and
the natural homomorphism $\mu_*\Cal O_Y(\ulcorner - F_t \urcorner) \to 
\mu_*\Cal O_E(\ulcorner - F_t \vert_E \urcorner)$
is surjective.
Since the negative part of $F_t$ is exceptional for $\mu$,
it follows that there is a natural isomorphism
$\Cal O_W \to \mu_*\Cal O_E(\ulcorner - F_t \vert_E \urcorner)$ if $t \ne 0$.

We may assume that there is a resolution of singularities 
$\sigma: V \to W$ which factors $\mu: E \to W$. 
Let $f: E \to V$ be the induced morphism:
$$
\CD
E @>>>     Y \\
@VfVV        \\
V            \\
@V{\sigma}VV  @VV{\mu}V\\
W @>>>     X.
\endCD
$$ 
We may also assume that there exist normal crossing divisors $P$ and $Q$ on 
$E$ and $V$, respectively, 
such that the conditions of Theorem 2 are satisfied 
for $F_t \vert_E$ if $t \ne 0$,
since we have
$$
\mu^*(K_X + D_t) \vert_E 
\sim_{\Bbb Q} (K_Y + E + F_t) \vert_E = K_E + F_t \vert_E. 
$$
We define $\Bbb Q$-divisors $M_t$ and $\Delta_t$ on $V$ for $0 \le t \le 1$
such that $K_E + F_t \vert_E \sim_{\Bbb Q} f^*(K_V + M_t + \Delta_t)$ as in
Theorem 2. 
By construction, the coefficients of $\Delta_t$ are less than $1$ for any $t$.
By Theorem 2, $M_t$ is nef for $t \ne 0$, hence for any $t$.

The surjectivity of the homomorphism
$\Cal O_W \to \mu_*\Cal O_E(\ulcorner - F_t \vert_E \urcorner)$ implies that,
if $\sigma_*Q_{\ell} \ne 0$, then 
there exists a $j$ such that $f(P_j) = Q_{\ell}$ and 
$d_j \ge 1 - w_{\ell j}$. 
Thus $0 \le \delta_{\ell}$ and $\sigma_*\Delta_t$ is effective.

We let $t = 0$, and set $M = M_0$ and $\Delta = \Delta_0$.
Since $M$ is nef, we may assume that there exists
rational numbers $q_{\ell}$ such that
$q_{\ell} > 0$ (resp. $= 0$) if $\sigma_*Q_{\ell} = 0$ (resp. $\ne 0$) and that
$M + \epsilon \sigma^*H - \epsilon' \sum_{\ell} q_{\ell}Q_{\ell}$ is ample for 
$0 < \epsilon' \ll \epsilon$. 
We take a general effective $\Bbb Q$-divisor
$M' \sim_{\Bbb Q} M + \epsilon \sigma^*H 
- \epsilon' \sum_{\ell} q_{\ell}Q_{\ell}$ with very small coefficients
and a very ample divisor as a support.
Let $D_W = \sigma_*(M' + \Delta)$.
Then we have $(K_X + D + \epsilon H) \vert_W \sim_{\Bbb Q} K_W + D_W$, and
$$
\sigma^*(K_W + D_W) 
\sim_{\Bbb Q} K_V + M' + \Delta + \epsilon' \sum_{\ell} q_{\ell}Q_{\ell}.
$$
If $\epsilon'$ is chosen small enough, 
then the coefficients on the right hand side are less than 1, and $(W, D_W)$ 
is KLT.
\hfill $\square$
\enddemo

\definition{Remark 3}
Since the choice of $M'$ is generic in the proof of Theorem 1, 
we can take $D_W$ such that the following holds.
Let $D'$ be an effective $\Bbb Q$-Cartier divisor on $X$ whose support 
does not contain $W$.
Assume that $(W, D_W + D' \vert_W)$ is not KLT.
Then $W$ is not a minimal center of log canonical singularities for 
$(X, D + D')$.
Indeed, the variation of Hodge structures considered in the proof does not 
change if we replace $D$ by $D + D'$. Then $\Delta$ is replaced by
$\Delta + \sigma^*(D' \vert_W)$, and $D_W$ by $D_W + D'\vert_W$.
\enddefinition

\definition{Remark 4}
There is an interesting application of the subadjunction formula to
the study on the existence of good ladders on Fano or log Fano varieties 
in [M] and [A].
Theorem 1 may also be related to Fujita's freeness conjecture.
\enddefinition

\enddocument